\documentstyle[12pt]{article}
\newcommand{\F}{{\cal F}}
\newcommand{\G}{{\cal G}}
\newcommand{\Oo}{{\cal O}}
\newcommand{\K}{{\cal K}}
\newcommand{\Hh}{{\cal H}}

\begin{document}
\begin{center}{\large {\bf Cohomology of coherent sheaves on a proper
scheme}}\bigskip\\
{\sc George R. Kempf\\
The Johns Hopkins University}\bigskip \bigskip \\
\end{center}

Let $f:X\rightarrow S$ be a proper morphism where $S = \mbox{Spec } A$ where
$A$ is a
noetherian ring.  Let $\F$ be a coherent sheaf on $X$.  Then a theorem of
Grothendieck
says that the cohomology groups $H^i(X,\F)$ are $A$-modules of finite type.  We
set out to
give an algorithm to compute these $A$-modules when $A$ is sufficiently
algorithmic.  We
found  general result about a $\Gamma(X,-)$-acylic resolution of $\F$.

We will first define a coherent sheaf $\G$ on $X$ to be primitive if support
$(\G)=\bar
x$ where $x$ is a point of $X$ and the natural mapping $\G \rightarrow
i_{x*}(\tilde \G_x)$
is injective where $i_x:\mbox{Spec } \Oo_{X,x}\rightarrow X$ is the natural
morphism.  The
associated points of $\G$ is denoted by $\mbox{Ass}(\G)$, which is a finite
set.

The theorem we prove is the following.\\

\noindent {\bf Theorem 1.}  {\it There exists an inclusion}
$\G\stackrel{\alpha}{\hookrightarrow}\oplus_{x \in \mbox{Ass}(\G)} \Hh_x$ {\it
where a) each
$\Hh_x$ is a primitive coherent sheaf with support $\bar x$, b) each $\Hh_x$ is
$\Gamma(X,-)$-acylic and its group of sections $\Gamma (X,{\cal H}_x)$ is an
$A$-module of finite type, and c) $\alpha$ is an isomorphism on an open dense
subset of the
support of $\G$.}

This theorem gives a method to compute the cohomology of $\F$.  One applies the
theorem
inductively to construct a $\Gamma(X,-)$-acylic coherent resolution
$0\rightarrow
\F\rightarrow \F^*$ and simply computes $H^i(X,\F)$ as the $i$-homology group
of the
complex $\Gamma(X, \F^*)$.

Thus Theorem 1 implies \\

\noindent {\bf Corollary} (Grothendieck).  {\it The
cohomology groups $H^i(X,\F)$ are $A$-modules of finite type.}\bigskip \\

\noindent {\bf \S1.  A lemma.}\\

We will use\\

\noindent {\bf Lemma 2.}  {\it If ${\cal F}$ is a coherent sheaf on $X$, there
exists an
embedding} $\F\stackrel{\beta}{\hookrightarrow}\oplus_{x \in \mbox{ Ass}(\F)}
\G_x$ {\it
where $\G_x$ is a primitive coherent sheaf with support $\bar x$ and $\beta$ is
an
isomorphism in an open dense subset of support $(\F)$.}\\

\noindent {\it Proof.}  We have an injection $\F \hookrightarrow \oplus_{x \in
\mbox{
Ass}({\cal F})} i_{x*}(\tilde {\F}_x)$ because the kernel would have no
associated points.
Similarly, we have the inclusion $\F_x \i\;
^{\mbox{\small{limit}}}_{\;\;\;\stackrel{\leftarrow}{n}} F_x / m^n_xF$ where
$m_x$ is the
maximal ideal of $\Oo_{x,x}$.  As $F$ is coherent it follows that
$\F\hookrightarrow
\oplus_{x \in \mbox{ Ass}(F)} i_{x*}(\widetilde{\F_x/m_x^{n_{x}}} \F_x)$ where
the $n_x$ are
sufficiently big.  This is our embedding if ${\cal G}_x$ is the image of ${\cal
F}$ in
$i_{x*}(\widetilde {F_x/m^{n_{x}}})$. \hfill QED\\

By the lemma to prove the theorem we may assume that $\G$ is primitive with
support $\bar x$
for some $x$ in $X$.\bigskip \\

\noindent {\bf \S2.  The proof.}\\

To prove the theorem we will use basically Grothendieck's idea of using Chow
coverings.
Let $Z$ be a closed subscheme of $X$ such that the support $(Z)=\bar x$ and
$\G$ admits the
structure of a $\Oo_Z$-module.  We will modify Chow lemma to the non-reduce
case.

Let $Z=\cup U_i$ be a finite affine covering of $Z$ by dense open subsets.  Let
$\bar
U_i\supset U_i$ be a projective imbedding for each $i$.  Take $W$ to be the
schematic
closure of $\cap U_i$ embedded in $X\times_S \prod \bar U_i$.  Then by the
usual reasoning the projection $W \rightarrow \Pi \bar U_i$ is
set-theoretically injective
proper morphism.  Thus $W$ is a projective scheme over $S$.

Let ${\cal M}$ be a very ample sheaf on $\Pi\bar U_i=K$.  Then ${\cal N} =
\pi^*_K{\cal
M}$ is very ample on $W$ and relatively ample for the projective morphism
$\pi^*_Z:W\rightarrow Z$, which is an isomorphism between open dense subsets.

Now let ${\cal L}_x$ be the maximal primitive quotient of $\pi^*_Z\G$ with
support $W$.  By
projective theory the $\K_x\equiv {\cal L}_x\otimes {\cal M}^{\otimes n}$ is
$\Gamma(W,
-)$-acylic and $\pi_{Z*}$-acylic for $n\gg 0$.  Thus $\pi_{Z*} (\K_*)=\Hh_x$ is
a
$\Gamma(X,-)$-acylic coherent primitive sheaf, by the degenerate Leray spectral
sequence.
It has support $\bar x$ and its sections are isomorphic to $\Gamma(W, \K_x)$,
which is an
$A$-module of finite type.

To get the inclusion $\F\rightarrow \Hh_x$ just multiply the natural map
$\F\rightarrow
\pi_{Z*} {\cal L}_x$ by a section of ${\cal M}^{\otimes n}$ which is general.
This
proves the theorem.\bigskip \\

\noindent {\bf References}
\begin{enumerate}

\item A. Grothendieck, \'El\'ements de G\'eom\'etrie Alg\'ebrique III
Publications
Mathematiques No.\ 111.

\item G. Kempf, Some elementary proofs of basis theorems in the cohomology of
quasi-coherent sheaves, Rocky Mountain Journal of Math., vol.\ 10, 1980, pp.\
637-645.
\end{enumerate}
\bigskip

\noindent {\small{\sc AMS (MOS) subject classification number 14F05.}}
\end{document}